\documentclass[11pt]{article}
\usepackage{amssymb,amsthm,amsmath}
\def\qed{\hfill $\Box$}

\newcommand\pf{\smallbreak\noindent \texttt{Proof}. }

\begin{document}

\newtheorem{thm}{Theorem}[section]
\newtheorem{prop}[thm]{Proposition}
\newtheorem{lem}[thm]{Lemma}
\newtheorem{cor}[thm]{Corollary}
\newtheorem{ex}[thm]{Example}
\renewcommand{\thefootnote}{*}

\title{\bf On some relationships between the center and the derived subalgebra in Poisson (2-3)-algebras}

\author{\textbf{P.~Ye.~Minaiev, O.~O.~Pypka}\\
Oles Honchar Dnipro National University, Dnipro, Ukraine\\
{\small e-mail: minaevp9595@gmail.com, sasha.pypka@gmail.com}\\
\textbf{I.~V.~Shyshenko}\\
Sumy State Pedagogical University named\\ after A.~S.~Makarenko, Sumy, Ukraine\\
{\small e-mail: shiinna@ukr.net}}
\date{}

\maketitle

\begin{abstract}
One of the classic results of group theory is the so-called Schur theorem. It states that if the central factor-group $G/\zeta(G)$ of a group $G$ is finite, then its derived subgroup $[G,G]$ is also finite. This result has numerous generalizations and modifications in group theory. At the same time, similar investigations were conducted in other algebraic structures, namely in modules, linear groups, topological groups, $n$-groups, associative algebras, Lie algebras, Lie $n$-algebras, Lie rings, Leibniz algebras. In 2021, L.A.~Kurdachenko, O.O.~Pypka and I.Ya.~Subbotin proved an analogue of Schur theorem for Poisson algebras: if the center of the Poisson algebra $P$ has finite codimension, then $P$ includes an ideal $K$ of finite dimension such that $P/K$ is abelian. In this paper, we continue similar studies for another algebraic structure. An analogue of Schur theorem for Poisson (2-3)-algebras is proved.
\end{abstract}



\thispagestyle{empty}

\section{Introduction.}
We start with a classical result of group theory. In 1951, B.H.~Neumann proved the so-called~\cite{KS2016} Schur theorem. This theorem states that if the central factor-group $G/\zeta(G)$ of a group $G$ is finite, then its derived subgroup $[G,G]$ is also finite~\cite{NB1951}. In other words, if $G/\zeta(G)$ is finite, then $G$ includes a finite normal subgroup $H=[G,G]$ such that $G/H$ is abelian. Moreover, there exists a function $w$ such that $|[G,G]|\leqslant w(t)$, where $t=|G/\zeta(G)|$. This theorem has numerous generalizations and modifications in group theory (see, for, example, \cite{DKP2014,DKP2015,KOP2015,KOP2016B,KP2014,P2017}). Furthermore, similar investigations were conducted in other algebraic structures of a different nature. Thus, analogues of Schur theorem were obtained for modules (see, for example, \cite{KSC2015}), linear groups~\cite{DKO2013}, topological groups~\cite{U1964}, $n$-groups~\cite{G2006}, associative algebras~\cite{RS2017}, Lie algebras~\cite{KPS2015,S1974}, Lie $n$-algebras~\cite{SV2014}, Lie rings~\cite{KPS2018}, Leibniz algebras~\cite{KOP2016A}.

A few years ago, L.A.~Kurdachenko, O.O.~Pypka and I.Ya.~Subbotin proved an analogue of Schur theorem for Poisson algebras~\cite{KPS2021}. More precisely, it was proved that if the center $\zeta(P)$ of the Poisson algebra $P$ has finite codimension, then $P$ includes an ideal $K$ of finite dimension such that $P/K$ is abelian. Moreover, they proved that $dim_{F}(K)\leqslant\frac{d(d^{2}-1)}{2}$.

Consider another algebraic structure that arises from Poisson algebras by replacing the binary operation of Lie bracket $[-,-]$ with the ternary operation $[-,-,-]$.

Let $P$ be a vector space over a field $F$. Then $P$ is called a \textit{Poisson $(2$-$3)$-algebra}, if $P$ has a binary operation $\cdot$ and ternary operation $[-,-,-]$ such that the product $\cdot$ forms a commutative associative algebra, the bracket $[-,-,-]$ forms a Lie $3$-algebra, and $[-,-,-]$ acts as a derivation of the product $\cdot$, that is $\cdot$ and $[-,-,-]$ satisfy the Leibniz ($2$-$3$)-identity. In other words,
\begin{gather*}
ab=ba, (ab)c=a(bc),\\
a(b+c)=ab+ac, (\lambda a)b=a(\lambda b)=\lambda(ab),\\
[a_{1}+a_{2},b,c]=[a_{1},b,c]+[a_{2},b,c],\\
[a,b_{1}+b_{2},c]=[a,b_{1},c]+[a,b_{2},c],\\
[a,b,c_{1}+c_{2}]=[a,b,c_{1}]+[a,b,c_{2}],\\
[\lambda a,b,c]=[a,\lambda b,c]=[a,b,\lambda c]=\lambda[a,b,c],\\
[a_{1},a_{2},a_{3}]=0 \mbox{ whenever }a_{i}=a_{j}\mbox{ for some }i\neq j, 1\leqslant i,j\leqslant 3,\\
[[a_{1},a_{2},a_{3}],b,c]=[[a_{1},b,c],a_{2},a_{3}]+[a_{1},[a_{2},b,c],a_{3}]+[a_{1},a_{2},[a_{3},b,c]],\\
[a_{1}a_{2},b,c]=a_{2}[a_{1},b,c]+a_{1}[a_{2},b,c]
\end{gather*}
for all $a,a_{1},a_{2},a_{3},b,b_{1},b_{2},c,c_{1},c_{2}\in P$, $\lambda\in F$. If we will consider $P$ as an associative and commutative algebra by multiplication $\cdot$, then we will denote it by $P(+,\cdot)$. If we will consider $P$ as a Lie 3-algebra by Lie 3-bracket $[-,-,-]$, then we will denote it by $P(+,[-,-,-])$.

This article is devoted to the extension of the statement, similar to Schur theorem, to Poisson (2-3)-algebras.

\section{Preliminary results.}
We begin with a series of important basic results that illustrate the properties of Poisson (2-3)-algebras.

\begin{lem}\label{L2.1}
Let $P$ be a Poisson $(2$-$3)$-algebra over a field $F$. Then for all elements of $P$
\begin{enumerate}
\item[\upshape(i)] $[a,b,c]=[b,c,a]=[c,a,b]=-[a,c,b]=-[b,a,c]=-[c,b,a]$;
\item[\upshape(ii)] $[a,b_{1}b_{2},c]=b_{2}[a,b_{1},c]+b_{1}[a,b_{2},c]$;
\item[\upshape(iii)] $[a,b,c_{1}c_{2}]=c_{2}[a,b,c_{1}]+c_{1}[a,b,c_{2}]$.
\end{enumerate}
\end{lem}
\pf (i) By definition, $[a_{1},a_{2},a_{3}]=0$ whenever $a_{i}=a_{j}$ for some $i\neq j$, $1\leqslant i,j\leqslant 3$. Since $[a,a,c]=0$ for all $a,c\in P$, then
\begin{gather*}
[a+b,a+b,c]=0,\\
[a,a,c]+[a,b,c]+[b,a,c]+[b,b,c]=0,\\
[a,b,c]=-[b,a,c].
\end{gather*}
Since $[a,b,a]=0$ for all $a,b\in P$, then
\begin{gather*}
[a+c,b,a+c]=0,\\
[a,b,a]+[a,b,c]+[c,b,a]+[c,b,c]=0,\\
[a,b,c]=-[c,b,a].
\end{gather*}
Since $[a,b,b]=0$ for all $a,b\in P$, then
\begin{gather*}
[a,b+c,b+c]=0,\\
[a,b,b]+[a,b,c]+[a,c,b]+[a,c,c]=0,\\
[a,b,c]=-[a,c,b].
\end{gather*}
Finally,
$$-[c,b,a]=-(-[b,c,a])=[b,c,a]$$
and
$$-[a,c,b]=-(-[c,a,b])=[c,a,b].$$

(ii) Since by (i) $[a,b,c]=[b,c,a]$, then
$$[a,b_{1}b_{2},c]=[b_{1}b_{2},c,a]=b_{2}[b_{1},c,a]+b_{1}[b_{2},c,a]=b_{2}[a,b_{1},c]+b_{1}[a,b_{2},c].$$

(iii) Since by (i) $[a,b,c]=[c,a,b]$, then
$$[a,b,c_{1}c_{2}]=[c_{1}c_{2},a,b]=c_{2}[c_{1},a,b]+c_{1}[c_{2},a,b]=c_{2}[a,b,c_{1}]+c_{1}[a,b,c_{2}].$$
\qed

\begin{lem}\label{L2.2}
Let $P$ be a Poisson $(2$-$3)$-algebra over a field $F$. Then for all elements of $P$
\begin{enumerate}
\item[\upshape(i)] $\left[\prod\limits_{1\leqslant i\leqslant k}a_{i},b,c\right]=\sum\limits_{1\leqslant i\leqslant k}\prod\limits_{\substack{1\leqslant\tau\leqslant k\\\tau\neq i}}a_{\tau}[a_{i},b,c]$;
\item[\upshape(ii)] $\left[a,\prod\limits_{1\leqslant j\leqslant s}b_{j},c\right]=\sum\limits_{1\leqslant j\leqslant s}\prod\limits_{\substack{1\leqslant\mu\leqslant s\\\mu\neq j}}b_{\mu}[a,b_{j},c]$;
\item[\upshape(iii)] $\left[a,b,\prod\limits_{1\leqslant r\leqslant t}c_{r}\right]=\sum\limits_{1\leqslant r\leqslant t}\prod\limits_{\substack{1\leqslant\nu\leqslant t\\\nu\neq r}}c_{\nu}[a,b,c_{r}]$;
\item[\upshape(iv)] $\left[\prod\limits_{1\leqslant i\leqslant k}a_{i},\prod\limits_{1\leqslant j\leqslant s}b_{j},\prod\limits_{1\leqslant r\leqslant t}c_{r}\right]=\sum\limits_{\substack{1\leqslant i\leqslant k\\1\leqslant j\leqslant s\\1\leqslant r\leqslant t}}\prod\limits_{\substack{1\leqslant\tau\leqslant k\\1\leqslant\mu\leqslant s\\1\leqslant\nu\leqslant t\\\tau\neq i,\mu\neq j,\nu\neq r}}a_{\tau}b_{\mu}c_{\nu}[a_{i},b_{j},c_{r}]$;
\item[\upshape(v)] $[a^{k},b^{s},c^{t}]=ksta^{k-1}b^{s-1}c^{t-1}[a,b,c]$.
\end{enumerate}
\end{lem}
\pf (i) Consider the product $\left[\prod\limits_{1\leqslant i\leqslant 3}a_{i},b,c\right]$:
\begin{gather*}
\left[\prod\limits_{1\leqslant i\leqslant 3}a_{i},b,c\right]=[a_{1}a_{2}a_{3},b,c]=[(a_{1}a_{2})a_{3},b,c]=\\
=a_{3}[a_{1}a_{2},b,c]+a_{1}a_{2}[a_{3},b,c]=a_{3}(a_{2}[a_{1},b,c]+a_{1}[a_{2},b,c])+a_{1}a_{2}[a_{3},b,c]=\\
=a_{2}a_{3}[a_{1},b,c]+a_{1}a_{3}[a_{2},b,c]+a_{1}a_{2}[a_{3},b,c]=\sum\limits_{1\leqslant i\leqslant 3}\prod\limits_{\substack{1\leqslant\tau\leqslant 3\\\tau\neq i}}a_{\tau}[a_{i},b,c].
\end{gather*}

Suppose we have already proved that
$$\left[\prod\limits_{1\leqslant i\leqslant d}a_{i},b,c\right]=\sum\limits_{1\leqslant i\leqslant d}\prod\limits_{\substack{1\leqslant\tau\leqslant d\\\tau\neq i}}a_{\tau}[a_{i},b,c].$$
Then
\begin{gather*}
\left[\prod\limits_{1\leqslant i\leqslant d+1}a_{i},b,c\right]=\left[\left(\prod\limits_{1\leqslant i\leqslant d}a_{i}\right)a_{d+1},b,c\right]=\\
=a_{d+1}\left[\prod\limits_{1\leqslant i\leqslant d}a_{i},b,c\right]+\prod\limits_{1\leqslant i\leqslant d}a_{i}[a_{d+1},b,c]=\\
=a_{d+1}\sum\limits_{1\leqslant i\leqslant d}\prod\limits_{\substack{1\leqslant\tau\leqslant d\\\tau\neq i}}a_{\tau}[a_{i},b,c]+\prod\limits_{1\leqslant i\leqslant d}a_{i}[a_{d+1},b,c]=\\
=\sum\limits_{1\leqslant i\leqslant d}\prod\limits_{\substack{1\leqslant\tau\leqslant d\\\tau\neq i}}(a_{\tau}a_{d+1})[a_{i},b,c]+\prod\limits_{1\leqslant i\leqslant d}a_{i}[a_{d+1},b,c]=\\
=\sum\limits_{1\leqslant i\leqslant d+1}\prod\limits_{\substack{1\leqslant\tau\leqslant d+1\\\tau\neq i}}a_{\tau}[a_{i},b,c].
\end{gather*}
Thus for every positive integer $k$ we proved that
$$\left[\prod\limits_{1\leqslant i\leqslant k}a_{i},b,c\right]=\sum\limits_{1\leqslant i\leqslant k}\prod\limits_{\substack{1\leqslant\tau\leqslant k\\\tau\neq i}}a_{\tau}[a_{i},b,c]$$

Similarly, using (ii) and (iii) of Lemma~\ref{L2.1}, we can prove that for all positive integers $s$ and $t$
\begin{gather*}
\textrm{(ii)} \left[a,\prod\limits_{1\leqslant j\leqslant s}b_{j},c\right]=\sum\limits_{1\leqslant j\leqslant s}\prod\limits_{\substack{1\leqslant\mu\leqslant s\\\mu\neq j}}b_{\mu}[a,b_{j},c];\\
\textrm{(iii)} \left[a,b,\prod\limits_{1\leqslant r\leqslant t}c_{r}\right]=\sum\limits_{1\leqslant r\leqslant t}\prod\limits_{\substack{1\leqslant\nu\leqslant t\\\nu\neq r}}c_{\nu}[a,b,c_{r}].
\end{gather*}

Combining (i), (ii) and (iii), we obtain that
\begin{gather*}
\textrm{(iv)} \left[\prod\limits_{1\leqslant i\leqslant k}a_{i},\prod\limits_{1\leqslant j\leqslant s}b_{j},\prod\limits_{1\leqslant r\leqslant t}c_{r}\right]=\sum\limits_{\substack{1\leqslant i\leqslant k\\1\leqslant j\leqslant s\\1\leqslant r\leqslant t}}\prod\limits_{\substack{1\leqslant\tau\leqslant k\\1\leqslant\mu\leqslant s\\1\leqslant\nu\leqslant t\\\tau\neq i,\mu\neq j,\nu\neq r}}a_{\tau}b_{\mu}c_{\nu}[a_{i},b_{j},c_{r}].
\end{gather*}

Finally, using (iv), we obtain that
\begin{gather*}
\textrm{(v)} [a^{k},b^{s},c^{t}]=ksta^{k-1}b^{s-1}c^{t-1}[a,b,c].
\end{gather*}
\qed

Let $P$ be a Poisson (2-3)-algebra over a field $F$. A subset $S$ of $P$ is called a \textit{subalgebra} of $P$ if $S$ is a subspace of $P$ and $ab,[a,b,c]\in S$ for all $a,b,c \in S$. A subset $I$ of $P$ is called an \textit{ideal} of $P$ if $I$ is a subspace of $P$ and $ab,[a,b,c]\in I$ for all $a\in I$ and $b,c\in P$.

We will say that a Poisson (2-3)-algebra $P$ is \textit{simple}, if it has only two ideals, $\langle0\rangle$ and $P$. As usual, $P$ is \textit{abelian}, if $[a,b,c]=0$ for all $a,b,c\in P$.

Let $P_{1}$ and $P_{2}$ be Poisson (2-3)-algebras over a field $F$. Then a mapping $f:P_{1}\rightarrow P_{2}$ is called a \textit{homomorphism}, if
\begin{gather*}
f(\lambda a)=\lambda f(a), f(a+b)=f(a)+f(b),\\
f(ab)=f(a)f(b), f([a,b,c])=[f(a),f(b),f(c)]
\end{gather*}
for all $a,b,c\in P_{1}$, $\lambda\in F$.

As usual, an injective homomorphism is called a \textit{monomorphism}, a surjective homomorphism is called an \textit{epimorphism}, and bijective homomorphism is called an \textit{isomorphism}.

\begin{prop}\label{P2.3}
Let $A$ be an arbitrary Poisson $(2$-$3)$-algebra over a field $F$. Then
\begin{enumerate}
\item[\upshape(i)] there exists a Poisson $(2$-$3)$-algebra $P$ over a field $F$ having a multiplicative identity element;
\item[\upshape(ii)] there exists a monomorphism $f:A\rightarrow P$ such that $\mathrm{Im}(f)$ is an ideal of $P$.
\end{enumerate}
\end{prop}
\pf (i) If $A$ has an identity element by multiplication $\cdot$, then all is proved. Therefore, suppose that $A$ has no an identity element by multiplication $\cdot$. Put $P=A\times F$ and define on $P$ the following operations:
\begin{gather*}
\lambda(a,\alpha)=(\lambda a,\lambda\alpha),\\
(a,\alpha)+(b,\beta)=(a+b,\alpha+\beta),\\
(a,\alpha)(b,\beta)=(ab+\alpha b+\beta a,\alpha\beta),\\
[(a,\alpha),(b,\beta),(c,\delta)]=([a,b,c],0_{F})
\end{gather*}
for all $a,b,c\in A$ and $\lambda,\alpha,\beta,\delta\in F$.

Let us show first that $P$ is an algebra. If $a,b,c\in A$ and $\alpha,\beta,\delta\in F$, then
\begin{gather*}
(a,\alpha)+(b,\beta)=(a+b,\alpha+\beta)=(b+a,\beta+\alpha)=(b,\beta)+(a,\alpha);
\end{gather*}
\begin{gather*}
((a,\alpha)+(b,\beta))+(c,\delta)=(a+b,\alpha+\beta)+(c,\delta)=\\
=((a+b)+c,(\alpha+\beta)+\delta)=(a+(b+c),\alpha+(\beta+\delta))=\\
=(a,\alpha)+(b+c,\beta+\delta)=(a,\alpha)+((b,\beta)+(c,\delta));
\end{gather*}
\begin{gather*}
(a,\alpha)+(0_{A},0_{F})=(a+0_{A},\alpha+0_{F})=(a,\alpha);
\end{gather*}
\begin{gather*}
(a,\alpha)+(-a,-\alpha)=(a+(-a),\alpha+(-\alpha))=(0_{A},0_{F});
\end{gather*}
\begin{gather*}
(a,\alpha)((b,\beta)+(c,\delta))=(a,\alpha)(b+c,\beta+\delta)=\\
=(a(b+c)+\alpha(b+c)+(\beta+\delta)a,\alpha(\beta+\delta))=\\
=((ab+\alpha b+\beta a)+(ac+\alpha c+\delta a),\alpha\beta+\alpha\delta)=\\
=(ab+\alpha b+\beta a,\alpha\beta)+(ac+\alpha c+\delta a,\alpha\delta)=\\
=(a,\alpha)(b,\beta)+(a,\alpha)(c,\delta);
\end{gather*}
\begin{gather*}
((a,\alpha)+(b,\beta))(c,\delta)=(a+b,\alpha+\beta)(c,\delta)=\\
=((a+b)c+(\alpha+\beta)c+\delta(a+b),(\alpha+\beta)\delta)=\\
=((ac+\alpha c+\delta a)+(bc+\beta c+\delta b),\alpha\delta+\beta\delta)=\\
=(ac+\alpha c+\delta a,\alpha\delta)+(bc+\beta c+\delta b,\beta\delta)=\\
=(a,\alpha)(c,\delta)+(b,\beta)(c,\delta).
\end{gather*}

Let $\lambda,\kappa\in F$. Then

\begin{gather*}
\lambda((a,\alpha)+(b,\beta))=\lambda(a+b,\alpha+\beta)=\\
=(\lambda(a+b),\lambda(\alpha+\beta))=(\lambda a+\lambda b,\lambda\alpha+\lambda\beta)=\\
=(\lambda a,\lambda\alpha)+(\lambda b,\lambda\beta)=\lambda(a,\alpha)+\lambda(b,\beta);
\end{gather*}
\begin{gather*}
(\lambda+\kappa)(a,\alpha)=((\lambda+\kappa)a,(\lambda+\kappa)\alpha)=\\
=(\lambda a+\kappa a,\lambda\alpha+\kappa\alpha)=(\lambda a,\lambda\alpha)+(\kappa a,\kappa\alpha)=\\
=\lambda(a,\alpha)+\kappa(a,\alpha);
\end{gather*}
\begin{gather*}
(\lambda\kappa)(a,\alpha)=((\lambda\kappa)a,(\lambda\kappa)\alpha)=\\
=(\lambda(\kappa a),\lambda(\kappa\alpha))=\lambda(\kappa a,\kappa\alpha)=\alpha(\kappa(a,\alpha));
\end{gather*}
\begin{gather*}
1_{F}(a,\alpha)=(1_{F}a,1_{F}\alpha)=(a,\alpha).
\end{gather*}

Finally,
\begin{gather*}
\lambda((a,\alpha)(b,\beta))=\lambda(ab+\alpha b+\beta a,\alpha\beta)=\\
=(\lambda ab+\lambda\alpha b+\lambda\beta a,\lambda\alpha\beta),\\
(\lambda(a,\alpha))(b,\beta)=(\lambda a,\lambda\alpha)(b,\beta)=\\
=(\lambda ab+\lambda\alpha b+\lambda\beta a,\lambda\alpha\beta),\\
(a,\alpha)(\lambda(b,\beta))=(a,\alpha)(\lambda b,\lambda\beta)=\\
=(\lambda ab+\lambda\alpha b+\lambda\beta a,\lambda\alpha\beta),
\end{gather*}
which implies that
$$\lambda((a,\alpha)(b,\beta))=(\lambda(a,\alpha))(b,\beta)=(a,\alpha)(\lambda(b,\beta)).$$

Therefore, $P$ is an algebra over a field $F$. Since
\begin{gather*}
(a,\alpha)(b,\beta)=(ab+\alpha b+\beta a,\alpha\beta)=\\
=(ba+\beta a +\alpha b,\beta\alpha)=(b,\beta)(a,\alpha),
\end{gather*}
$P$ is commutative. Moreover,
\begin{gather*}
((a,\alpha)(b,\beta))(c,\delta)=(ab+\alpha b+\beta a,\alpha\beta)(c,\delta)=\\
=(abc+\alpha bc +\beta ac+\alpha\beta c+\delta ab+\alpha\delta b+\beta\delta a,\alpha\beta\delta)
\end{gather*}
and
\begin{gather*}
(a,\alpha)((b,\beta)(c,\delta))=(a,\alpha)(bc+\beta c+\delta b,\beta\delta)=\\
=(abc+\beta ac+\delta ab+\alpha bc +\alpha\beta c+\alpha\delta b+\beta\delta a,\alpha\beta\delta),
\end{gather*}
which implies that
\begin{gather*}
((a,\alpha)(b,\beta))(c,\delta)=(a,\alpha)((b,\beta)(c,\delta)).
\end{gather*}
Thus, $P$ is an associative algebra. Since
$$[(a_{1},\alpha_{1}),(a_{2},\alpha_{2}),(a_{3},\alpha_{3})]=([a_{1},a_{2},a_{3}],0_{F}),$$
then
$$[(a_{1},\alpha_{1}),(a_{2},\alpha_{2}),(a_{3},\alpha_{3})]=(0_{A},0_{F})=0_{P}$$
whenever $(a_{i},\alpha_{i})=(a_{j},\alpha_{j})$ for some $i\neq j$, $1\leqslant i,j\leqslant 3$.

Furthermore,
\begin{gather*}
[[(a_{1},\alpha_{1}),(a_{2},\alpha_{2}),(a_{3},\alpha_{3})],(b,\beta),(c,\delta)]=\\
=[([a_{1},a_{2},a_{3}],0_{F}),(b,\beta),(c,\delta)]=([[a_{1},a_{2},a_{3}],b,c],0_{F}),
\end{gather*}
and
\begin{gather*}
[[(a_{1},\alpha_{1}),(b,\beta),(c,\delta)],(a_{2},\alpha_{2}),(a_{3},\alpha_{3})]+\\
+[(a_{1},\alpha_{1}),[(a_{2},\alpha_{2}),(b,\beta),(c,\delta)],(a_{3},\alpha_{3})]+\\
+[(a_{1},\alpha_{1}),(a_{2},\alpha_{2}),[(a_{3},\alpha_{3}),(b,\beta),(c,\delta)]]=\\
=[([a_{1},b,c],0_{F}),(a_{2},\alpha_{2}),(a_{3},\alpha_{3})]+\\
+[(a_{1},\alpha_{1}),([a_{2},b,c],0_{F}),(a_{3},\alpha_{3})]+\\
+[(a_{1},\alpha_{1}),(a_{2},\alpha_{2}),([a_{3},b,c],0_{F})]=\\
=([[a_{1},b,c],a_{2},a_{3}],0_{F})+([a_{1},[a_{2},b,c],a_{3}],0_{F})+([a_{1},a_{2},[a_{3},b,c]],0_{F})=\\
=([[a_{1},b,c],a_{2},a_{3}]+[a_{1},[a_{2},b,c],a_{3}]+[a_{1},a_{2},[a_{3},b,c]],0_{F}),
\end{gather*}
which shows that
\begin{gather*}
[[(a_{1},\alpha_{1}),(a_{2},\alpha_{2}),(a_{3},\alpha_{3})],(b,\beta),(c,\delta)]=\\
=[[(a_{1},\alpha_{1}),(b,\beta),(c,\delta)],(a_{2},\alpha_{2}),(a_{3},\alpha_{3})]+\\
+[(a_{1},\alpha_{1}),[(a_{2},\alpha_{2}),(b,\beta),(c,\delta)],(a_{3},\alpha_{3})]+\\
+[(a_{1},\alpha_{1}),(a_{2},\alpha_{2}),[(a_{3},\alpha_{3}),(b,\beta),(c,\delta)]].
\end{gather*}

Consider the Leibniz (2-3)-identity for $P$.
\begin{gather*}
[(a_{1},\alpha_{1})(a_{2},\alpha_{2}),(b,\beta),(c,\delta)]=\\
=[(a_{1}a_{2}+\alpha_{1}a_{2}+\alpha_{2}a_{1},\alpha_{1}\alpha_{2}),(b,\beta),(c,\delta)]=\\
=([a_{1}a_{2}+\alpha_{1}a_{2}+\alpha_{2}a_{1},b,c],0_{F})=\\
=([a_{1}a_{2},b,c]+\alpha_{1}[a_{2},b,c]+\alpha_{2}[a_{1},b,c],0_{F}).
\end{gather*}
On the other hand,
\begin{gather*}
(a_{2},\alpha_{2})[(a_{1},\alpha_{1}),(b,\beta),(c,\delta)]+(a_{1},\alpha_{1})[(a_{2},\alpha_{2}),(b,\beta),(c,\delta)]=\\
=(a_{2},\alpha_{2})([a_{1},b,c],0_{F})+(a_{1},\alpha_{1})([a_{2},b,c],0_{F})=\\
=(a_{2}[a_{1},b,c]+\alpha_{2}[a_{1},b,c]+0_{F}a_{2},\alpha_{2}0_{F})+(a_{1}[a_{2},b,c]\\
+\alpha_{1}[a_{2},b,c]+0_{F}a_{1},\alpha_{1}0_{F})=\\
=(a_{2}[a_{1},b,c]+\alpha_{2}[a_{1},b,c],0_{F})+(a_{1}[a_{2},b,c]+\alpha_{1}[a_{2},b,c],0_{F})=\\
=(a_{2}[a_{1},b,c]+a_{1}[a_{2},b,c]+\alpha_{1}[a_{2},b,c]+\alpha_{2}[a_{1},b,c],0_{F}).
\end{gather*}
Since $[a_{1}a_{2},b,c]=a_{2}[a_{1},b,c]+a_{1}[a_{2},b,c]$,
\begin{gather*}
[(a_{1},\alpha_{1})(a_{2},\alpha_{2}),(b,\beta),(c,\delta)]=\\
=(a_{2},\alpha_{2})[(a_{1},\alpha_{1}),(b,\beta),(c,\delta)]+(a_{1},\alpha_{1})[(a_{2},\alpha_{2}),(b,\beta),(c,\delta)].
\end{gather*}

Taking into account the above, we can say that $P$ is a Poisson (2-3)-algebra.

Finally, $(0_{A},1_{F})$ is an identity element by multiplication $\cdot$:
\begin{gather*}
(a,\alpha)(0_{A},1_{F})=(a0_{A}+\alpha0_{A}+1_{F}a,\alpha1_{F})=(a,\alpha).
\end{gather*}

(ii) Consider the mapping $f:A\rightarrow P$, which defined by the rule $f(a)=(a,0_{F})$ for all $a\in A$. Then
\begin{gather*}
f(\lambda a)=(\lambda a,0_{F})=(\lambda a,\lambda0_{F})=\lambda(a,0_{F})=\lambda f(a);\\
f(a+b)=(a+b,0_{F})=(a+b,0_{F}+0_{F})=(a,0_{F})+(b,0_{F})=f(a)+f(b);\\
f(ab)=(ab,0_{F})=(ab+0_{F}b+0_{F}a,0_{F}0_{F})=(a,0_{F})(b,0_{F})=f(a)f(b);\\
f([a,b,c])=([a,b,c],0_{F})=[(a,0_{F}),(b,0_{F}),(c,0_{F})]=[f(a),f(b),f(c)].
\end{gather*}
Thus $f$ is a homomorphism. Clearly $f$ is injective, so that $f$ is a monomorphism. Let
$$(a,0_{F}),(b,0_{F})\in\mathrm{Im}(f)=\{(x,0_{F})|\ x\in A\}.$$
Since
\begin{gather*}
(a,0_{F})-(b,0_{F})=(a-b,0_{F})\in\mathrm{Im}(f),\\
\lambda(a,0_{F})=(\lambda a,0_{F})\in\mathrm{Im}(f),
\end{gather*}
$\mathrm{Im}(f)$ is a subspace of $P$. Finally, if $(a,0_{F})\in\mathrm{Im}(f)$, $(b,\beta),(c,\delta)\in P$, then
\begin{gather*}
(a,0_{F})(b,\beta)=(ab+0_{F}b+\beta a,0_{F}\beta)=(ab+\beta a,0_{F})\in\mathrm{Im}(f),\\
[(a,0_{F}),(b,\beta),(c,\delta)]=([a,b,c],0_{F})\in\mathrm{Im}(f),
\end{gather*}
which shows that $\mathrm{Im}(f)$ is an ideal of $P$.
\qed

Proposition~\ref{P2.3} shows that we can consider only Poisson (2-3)-algebras with an identity element $1_{P}$ by multiplication $\cdot$.

Let $A,B,C$ be subspaces of a Poisson (2-3)-algebra $P$. Denote by
\begin{itemize}
\item $A+B$ the subspace of $P$ with elements of the form $a+b$, $a\in A$, $b\in B$;
\item $AB$ the subspace of $P$, which generated by the subset $\{ab|\ a\in A,b\in B\}$;
\item $[A,B,C]$ the subspace of $P$, which generated by the subset $\{[a,b,c]|\ a\in A,b\in B,c\in C\}$.
\end{itemize}

Clearly, every element of $AB$ has a following form
$$a_{1}b_{1}+\ldots+a_{n}b_{n}$$
where $a_{1},\ldots a_{n}\in A$, $b_{1},\ldots b_{n}\in B$.

Similarly, every element of $[A,B,C]$ has a following form
$$[a_{1},b_{1},c_{1}]+\ldots+[a_{n},b_{n},c_{n}]$$
where $a_{1},\ldots a_{n}\in A$, $b_{1},\ldots b_{n}\in B$ and $c_{1},\ldots c_{n}\in C$.

\begin{prop}\label{P2.4}
Let $P$ be a Poisson $(2$-$3)$-algebra over a field $F$.
\begin{enumerate}
\item[\upshape(i)] If $A$ is a subalgebra of $P$ and $B$ is an ideal of $P$, then $A+B$ is a subalgebra of $P$.
\item[\upshape(ii)] If $A$ is a subalgebra of $P$ and $B$ is an ideal of $P$, then $AB$ is a subalgebra of $P$.
\item[\upshape(iii)] If $A$, $B$ are ideals of $P$, then $A+B$ is an ideal of $P$.
\item[\upshape(iv)] If $A$, $B$ are ideals of $P$, then $AB$ is an ideal of $P$.
\item[\upshape(v)] If $P$ is non-simple, then $P$ has a proper non-zero maximal ideal.
\item[\upshape(vi)] If $a\in P$, then $aP=\{ab|\ b\in P\}$ is a subalgebra of $P$ and an ideal of $P(+,\cdot)$.
\end{enumerate}
\end{prop}
\pf (i) If $x,y,z\in A+B$, then $x=a_{1}+b_{1}$, $y=a_{2}+b_{2}$, $z=a_{3}+b_{3}$ where $a_{1},a_{2},a_{3}\in A$, $b_{1},b_{2},b_{3}\in B$. If $\lambda\in F$, then
\begin{gather*}
\lambda x=\lambda(a_{1}+b_{1})=\lambda a_{1}+\lambda b_{1}\in A+B;\\
x-y=(a_{1}+b_{1})-(a_{2}+b_{2})=(a_{1}-a_{2})+(b_{1}-b_{2})\in A+B;\\
xy=(a_{1}+b_{1})(a_{2}+b_{2})=a_{1}a_{2}+(a_{1}b_{2}+b_{1}a_{2}+b_{1}b_{2})\in A+B;\\
[x,y,z]=[a_{1}+b_{1},a_{2}+b_{2},a_{3}+b_{3}]=\\
=[a_{1},a_{2},a_{3}]+([a_{1},a_{2},b_{3}]+[a_{1},b_{2},a_{3}]+[a_{1},b_{2},b_{3}]+\\
+[b_{1},a_{2},a_{3}]+[b_{1},a_{2},b_{3}]+[b_{1},b_{2},a_{3}]+[b_{1},b_{2},b_{3}])\in A+B.
\end{gather*}
Thus $A+B$ is a subalgebra of $P$.

(ii) If $x,y,z\in AB$, then
\begin{gather*}
x=a_{1}b_{1}+\ldots+a_{n}b_{n}=\sum\limits_{1\leqslant i\leqslant n}a_{i}b_{i},\\
y=c_{1}d_{1}+\ldots+c_{s}d_{s}=\sum\limits_{1\leqslant j\leqslant s}c_{j}d_{j},\\
z=h_{1}k_{1}+\ldots+h_{t}k_{t}=\sum\limits_{1\leqslant r\leqslant t}h_{r}k_{r}
\end{gather*}
where $a_{i},c_{j},h_{r}\in A$, $b_{i},d_{j},k_{r}\in B$, $1\leqslant i\leqslant n$, $1\leqslant j\leqslant s$, $1\leqslant r\leqslant t$. If $\lambda\in F$, then
\begin{gather*}
\lambda x=\lambda(a_{1}b_{1}+\ldots+a_{n}b_{n})=(\lambda a_{1})b_{1}+\ldots+(\lambda a_{n})b_{n}\in AB;\\
x-y=(a_{1}b_{1}+\ldots+a_{n}b_{n})-(c_{1}d_{1}+\ldots+c_{s}d_{s})\in AB;\\
xy=\left(\sum\limits_{1\leqslant i\leqslant n}a_{i}b_{i}\right)\left(\sum\limits_{1\leqslant j\leqslant s}c_{j}d_{j}\right)=\sum\limits_{\substack{1\leqslant i\leqslant n\\1\leqslant j\leqslant s}}(a_{i}b_{i})(c_{j}d_{j})\\
=\sum\limits_{\substack{1\leqslant i\leqslant n\\1\leqslant j\leqslant s}}(a_{i}c_{j})(b_{i}d_{j})\in AB.
\end{gather*}

Consider the product $[x,y,z]$:
$$[x,y,z]=\left[\sum\limits_{1\leqslant i\leqslant n}a_{i}b_{i},\sum\limits_{1\leqslant j\leqslant s}c_{j}d_{j},\sum\limits_{1\leqslant r\leqslant t}h_{r}k_{r}\right]=\sum\limits_{\substack{1\leqslant i\leqslant n\\1\leqslant j\leqslant s\\1\leqslant r\leqslant t}}[a_{i}b_{i},c_{j}d_{j},h_{r}k_{r}].$$
Taking into account Lemma~\ref{L2.2}~(iv), we have
\begin{gather*}
[a_{i}b_{i},c_{j}d_{j},h_{r}k_{r}]=\\
=b_{i}d_{j}k_{r}[a_{i},c_{j},h_{r}]+b_{i}d_{j}h_{r}[a_{i},c_{j},k_{r}]+b_{i}c_{j}k_{r}[a_{i},d_{j},h_{r}]+b_{i}c_{j}h_{r}[a_{i},d_{j},k_{r}]+\\
+a_{i}d_{j}k_{r}[b_{i},c_{j},h_{r}]+a_{i}d_{j}h_{r}[b_{i},c_{j},k_{r}]+a_{i}c_{j}k_{r}[b_{i},d_{j},h_{r}]+a_{i}c_{j}h_{r}[b_{i},d_{j},k_{r}].
\end{gather*}
Since $A$ is a subalgebra of $P$, then $[a_{i},c_{j},h_{r}],c_{j}h_{r},a_{i}h_{r},a_{i}c_{j},a_{i}c_{j}h_{r}\in A$. On the other hand, since $B$ is an ideal of $P$, then $b_{i}d_{j}k_{r}$, $b_{i}d_{j}[a_{i},c_{j},k_{r}]$, $b_{i}k_{r}[a_{i},d_{j},h_{r}]$, $b_{i}[a_{i},d_{j},k_{r}]$, $d_{j}k_{r}[b_{i},c_{j},h_{r}]$, $d_{j}[b_{i},c_{j},k_{r}]$, $k_{r}[b_{i},d_{j},h_{r}]$, $[b_{i},d_{j},k_{r}]\in B$. Thus each term in the expansion of $[a_{i}b_{i},c_{j}d_{j},h_{r}k_{r}]$ belongs to $AB$. Obviously, this is true for all $i,j,r$, $1\leqslant i\leqslant n$, $1\leqslant j\leqslant s$, $1\leqslant r\leqslant t$, which implies that $[x,y,z]\in AB$. Therefore, $AB$ is a subalgebra of $P$.

(iii) As above we can show that $A+B$ is a subspace of $P$. If $x\in A+B$, then $x=a+b$, $a\in A$, $b\in B$. Let $y,z\in P$. Then
\begin{gather*}
xy=(a+b)y=ay+by\in A+B;\\
[x,y,z]=[a+b,y,z]=[a,y,z]+[b,y,z]\in A+B.
\end{gather*}
Thus $A+B$ is an ideal of $P$.

(iv) As above we can show that $AB$ is a subspace of $P$. If $x\in AB$, then
$$x=a_{1}b_{1}+\ldots+a_{n}b_{n}=\sum\limits_{1\leqslant i\leqslant n}a_{i}b_{i},$$
$a_{i}\in A$, $b_{i}\in B$, $1\leqslant i\leqslant n$. Let $y,z\in P$. Then
\begin{gather*}
xy=\left(\sum\limits_{1\leqslant i\leqslant n}a_{i}b_{i}\right)y=\sum\limits_{1\leqslant i\leqslant n}a_{i}(b_{i}y)\in AB;\\
[x,y,z]=\left[\sum\limits_{1\leqslant i\leqslant n}a_{i}b_{i},y,z\right]=\sum\limits_{1\leqslant i\leqslant n}[a_{i}b_{i},y,z]\\
=\sum\limits_{1\leqslant i\leqslant n}(b_{i}[a_{i},y,z]+a_{i}[b_{i},y,z])\in AB.
\end{gather*}
Therefore, $AB$ is an ideal of $P$.

(v) If $P$ is not simple, then $P$ includes a proper non-zero ideal. Put
$$\mathfrak{S}=\{A|\ A\ \mbox{is a proper non-zero ideal of }P\}.$$
Obviously, $\mathfrak{S}$ is non-empty. Since $1_{P}\not\in A$ for each $A\in\mathfrak{S}$, an union of every linearly ordered (by inclusion) subset of $\mathfrak{S}$ belongs to $\mathfrak{S}$. By Zorn's Lemma, a family $\mathfrak{S}$ has a maximal element.

(vi) If $x,y,z\in aP=\{ab|\ b\in P\}$, then $x=ab_{1}$, $y=ab_{2}$ and $z=ab_{3}$. Let $\lambda\in F$, then
\begin{gather*}
\lambda x=\lambda ab_{1}=a(\lambda b_{1})\in aP;\\
x-y=ab_{1}-ab_{2}=a(b_{1}-b_{2})\in aP;\\
xy=(ab_{1})(ab_{2})=a(b_{1}ab_{2})\in aP.
\end{gather*}
In particular, since for any $c\in P$ we have $xc=(ab_{1})c=a(b_{1}c)\in P$, then $aP$ is an ideal of $P(+,\cdot)$.

Finally, using Lemma~\ref{L2.2}~(iv), we have
\begin{gather*}
[ab_{1},ab_{2},ab_{3}]=\\
=b_{1}b_{2}b_{3}[a,a,a]+b_{1}b_{2}a[a,a,b_{3}]+b_{1}ab_{3}[a,b_{2},a]+b_{1}aa[a,b_{2},b_{3}]+\\
+ab_{2}b_{3}[b_{1},a,a]+ab_{2}a[b_{1},a,b_{3}]+aab_{3}[b_{1},b_{2},a]+aaa[b_{1},b_{2},b_{3}]=\\
=b_{1}aa[a,b_{2},b_{3}]+ab_{2}a[b_{1},a,b_{3}]+aab_{3}[b_{1},b_{2},a]+aaa[b_{1},b_{2},b_{3}]=\\
=a(ab_{1}[a,b_{2},b_{3}]+ab_{2}[b_{1},a,b_{3}]+ab_{3}[b_{1},b_{2},a]+a^{2}[b_{1},b_{2},b_{3}])\in aP.
\end{gather*}
The last inclusion shows that $aP$ is a subalgebra of $P$ for each $a\in P$.
\qed

Let $P$ be a Poisson (2-3)-algebra over a field $F$. Put
$$\zeta(P)=\{a\in P|\ [a,b,c]=0\ \mbox{for all }b,c\in P\}.$$
The subset $\zeta(P)$ is called the \textit{center} of $P$. We note that $\zeta(P)$ is an ideal of Lie 3-algebra $P(+,[-,-,-])$.

\begin{prop}\label{P2.5}
Let $P$ be a Poisson $(2$-$3)$-algebra over a field $F$. Then
\begin{enumerate}
\item[\upshape(i)] $\zeta(P)$ is a subalgebra of $P$;
\item[\upshape(ii)] $\zeta(P)$ contains every idempotent of $P$, in particular, $1_{P}\in\zeta(P)$.
\end{enumerate}
\end{prop}
\pf (i) Let $\lambda\in F$, $a\in\zeta(P)$, $b,c\in P$. Then
$$[\lambda a,b,c]=\lambda[a,b,c]=0.$$

If $a_{1},a_{2}\in\zeta(P)$ and $b,c,\in P$, then
$$[a_{1}-a_{2},b,c]=[a_{1},b,c]-[a_{2},b,c]=0,$$
which shows that $a_{1}-a_{2}\in\zeta(P)$, so that $\zeta(P)$ is a subspace of $P$.

Moreover,
$$[a_{1}a_{2},b,c]=a_{2}[a_{1},b,c]+a_{1}[a_{2},b,c]=0,$$
so that $a_{1}a_{2}\in\zeta(P)$.

Finally, let $a_{1},a_{2},a_{3}\in\zeta(P)$ and $b,c,\in P$. Then
\begin{gather*}
[[a_{1},a_{2},a_{3}],b,c]=[[a_{1},b,c],a_{2},a_{3}]+[a_{1},[a_{2},b,c],a_{3}]+[a_{1},a_{2},[a_{3},b,c]]=\\
=[0,a_{2},a_{3}]+[a_{1},0,a_{3}]+[a_{1},a_{2},0]=0,
\end{gather*}
which implies that $[a_{1},a_{2},a_{3}]\in\zeta(P)$. Thus $\zeta(P)$ is a subalgebra of $P$.

(ii) Let $e$ be an arbitrary idempotent of $P$. Let $b,c\in P$. Using Lemma~\ref{L2.2}~(v) we have
$$[e,b,c]=[e^{2},b,c]=2e[e,b,c].$$
If $2e[e,b,c]=0$, then all is proved. Suppose that $2e[e,b,c]\neq0$. Then
$$2e[e,b,c]-[e,b,c]=0.$$
It follows that
\begin{gather*}
0=2e(2e[e,b,c]-[e,b,c])=4e^{2}[e,b,c]-2e[e,b,c]=\\
=4e[e,b,c]-2e[e,b,c]=2e[e,b,c],
\end{gather*}
and we obtain a contradiction. Thus $2e[e,b,c]=0$, so that $e\in\zeta(P)$. In particular, since $1_{P}=1_{P}^{2}$, $1_{P}\in\zeta(P)$.
\qed

\section{Analogue of Schur theorem for Poisson (2-3)-algebras}
\begin{lem}\label{L3.1}
Let $L$ be a Lie 3-algebra over a field $F$. Suppose that the factor-algebra $L/\zeta(L)$ has a finite dimension $d$ and let $\{e_{1}+\zeta(L),\ldots,e_{d}+\zeta(L)\}$ be a basis of $L/\zeta(L)$. Then $[L,L,L]$ generated by the elements $[u_{1},u_{2},u_{3}]$ where $u_{j}\in\{e_{1},\ldots,e_{d}\}$, $1\leqslant j\leqslant 3$.
\end{lem}
\pf Put $Z=\zeta(L)$. Then for every $x\in L$ we have $x=\lambda_{1}e_{1}+\ldots+\lambda_{d}e_{d}+s_{x}$ for some $\lambda_{1},\ldots,\lambda_{d}\in F$ and $s_{x}\in Z$. Let $y=\mu_{1}e_{1}+\ldots+\mu_{d}e_{d}+s_{y}$ and $z=\nu_{1}e_{1}+\ldots+\nu_{d}e_{d}+s_{z}$, $\mu_{1},\ldots,\mu_{d},\nu_{1},\ldots,\nu_{d}\in F$, $s_{y},s_{z}\in Z$. Then
\begin{gather*}
[x,y,z]\\
=[\lambda_{1}e_{1}+\ldots+\lambda_{d}e_{d}+s_{x},\mu_{1}e_{1}+\ldots+\mu_{d}e_{d}+s_{y},\nu_{1}e_{1}+\ldots+\nu_{d}e_{d}+s_{z}]\\
=\sum\limits_{1\leqslant i,j,k\leqslant d}\lambda_{i}\mu_{j}\nu_{k}[e_{i},e_{j},e_{k}].
\end{gather*}
As we can see, $[L,L,L]$ is a subalgebra generates by the elements $[e_{i},e_{j},e_{k}]$, $1\leqslant i,j,k\leqslant d$. More precisely, since $[e_{i},e_{j},e_{k}]=sign(\sigma)[e_{\sigma(i)},e_{\sigma(j)},e_{\sigma(k)}]$ ($\sigma\in S_{3}$) and $[e_{i},e_{j},e_{k}]=0$ whenever $e_{m}=e_{n}$ for some $m\neq n$, $1\leqslant i,j,k\leqslant d$, $[L,L,L]$ generates by the elements $[e_{i},e_{j},e_{k}]$ where $1\leqslant i<j<k\leqslant d$.
\qed

\begin{thm}\label{T3.2}
Let $P$ be a Poisson $(2$-$3)$-algebra over a field $F$. Suppose that the center of $P$ has a finite codimension $d$. Then $P$ includes an ideal $K$ of finite dimension at most $\frac{d(d^{2}-1)(d-2)}{6}$ such that $P/K$ is abelian.
\end{thm}
\pf Put $Z=\zeta(P)$. Then $P=Z\oplus A$ for some subspace $A$ of $P$. Choose a basis $\{e_{1},\ldots,e_{d}\}$ in the subspace $A$. Then for every element $x\in P$ we have
$$x=\lambda_{1}e_{1}+\ldots+\lambda_{d}e_{d}+z_{x}$$
for some $\lambda_{1},\ldots,\lambda_{d}\in F$ and $z_{x}\in Z$.

A subspace $[P,P,P]$ is an ideal of a Lie 3-algebra $P(+,[-,-,-])$. Lemma~\ref{L3.1} shows that $[P,P,P]$ generates as a subspace by the elements $[e_{i},e_{j},e_{k}]$ where $1\leqslant i<j<k\leqslant d$.

Consider an ideal $K$ of an associative algebra $P(+,\cdot)$ generated by $[P,P,P]$. Every element from $K$ has a form $a_{1}x_{1}+\ldots+a_{r}x_{r}$ where $a_{1},\ldots,a_{r}\in[P,P,P]$, $x_{1},\ldots,x_{r}$ are the arbitrary elements of $P$. Let $x$ be an arbitrary element of $P$, $x=\lambda_{1}e_{1}+\ldots+\lambda_{d}e_{d}+z_{x}$ where $\lambda_{1},\ldots,\lambda_{d}\in F$ and $z_{x}\in Z$. We have
\begin{gather*}
[e_{i},e_{j},e_{k}]x=[e_{i},e_{j},e_{k}](\lambda_{1}e_{1}+\ldots+\lambda_{d}e_{d}+z_{x})\\
=\lambda_{1}e_{1}[e_{i},e_{j},e_{k}]+\ldots+\lambda_{d}e_{d}[e_{i},e_{j},e_{k}]+z_{x}[e_{i},e_{j},e_{k}].
\end{gather*}
Using Leibniz rule we obtain
$$[e_{i},e_{j},e_{k}z_{x}]=z_{x}[e_{i},e_{j},e_{k}]+e_{k}[e_{i},e_{j},z_{x}]=z_{x}[e_{i},e_{j},e_{k}].$$
For element $e_{k}z_{x}$ we have the decomposition $e_{k}z_{x}=\nu_{1}e_{1}+\ldots+\nu_{d}e_{d}+z_{x,k}$. Therefore
$$[e_{i},e_{j},e_{k}z_{x}]=[e_{i},e_{j},\nu_{1}e_{1}+\ldots+\nu_{d}e_{d}+z_{x,k}]=\nu_{1}[e_{i},e_{j},e_{1}]+\ldots+\nu_{d}[e_{i},e_{j},e_{d}].$$

These equalities show that $K$ as a vector space is generated by the elements $[e_{i},e_{j},e_{k}]$, $e_{s}[e_{i},e_{j},e_{k}]$, $1\leqslant s\leqslant d$, $1\leqslant i<j<k\leqslant d$. It follows that $K$ has a dimension at most
\begin{gather*}
\frac{d(d-1)(d-2)}{6}+d\frac{d(d-1)(d-2)}{6}=\frac{d(d-1)(d-2)}{6}(d+1)\\
=\frac{d(d^{2}-1)(d-2)}{6}.
\end{gather*}
The inclusion $[K,K,K]\leqslant[P,P,P]\leqslant K$ shows that $K$ is a subalgebra of $P$.

Moreover, $[K,P,P]\leqslant[P,P,P]\leqslant K$, so that $K$ is an ideal of Lie 3-algebra $P(+,[-,-,-])$.
\qed

\end{document}